\documentclass[12pt]{article}
\usepackage{geometry}                
\usepackage{tikz}
\usepackage{graphicx}
\usepackage{amssymb}					
\usepackage{amsmath}					
\usepackage{mathrsfs}
\usepackage{epstopdf}
\usepackage{amsthm}
\usepackage{float}

\def\b0{{\bf 0}}
\def\b1{{\bf 1}}

\def\cC{{\cal C}}

\def\cF{{\cal F}}
\def\cH{{\cal H}}

\def\cC{{\cal C}}

\def\cH{{\cal H}}

\def\n{\noindent}

\usepackage{float}
\usepackage{enumitem}
\definecolor{forestgreen}{rgb}{0,0.5,0}
\definecolor{silver}{rgb}{0.85,0.85,0.85}
\definecolor{lightsilver}{rgb}{0.97,0.97,0.97}
\definecolor{rust}{rgb}{0.75,0,0}
\definecolor{lightpink}{RGB}{253,231,232}

\newtheorem*{theorem*}{Theorem}
\newtheorem*{lemma*}{Lemma}
\newtheorem*{example*}{Example}
\newtheorem{theorem}{Theorem}
\newtheorem{lemma}{Lemma}

\newtheorem{proposition}{Proposition}

\usepackage{marginnote}

\begin{document}
\title{
Euler's Theorem for Regular CW-Complexes}
\author{Richard H. Hammack\footnote{Supported by Simons Collaboration Grant for Mathematicians
523748.}\\
\texttt{rhammack@vcu.edu}\\
\\
Paul C. Kainen\footnote{Corresponding author $\;\;\;$ this article will appear in {\it Combinatorica}, 2024}\\
 \texttt{kainen@georgetown.edu}
 }

\newcommand{\Addresses}{{
  \bigskip
  \footnotesize

\noindent
Paul C. Kainen, \textsc{Department of Mathematics and Statistics,\\ Georgetown University, 
37th and O Streets, N.W., Washington DC 20057 USA}\\
\vspace{-0.27cm}

\medskip
\noindent
Richard H. Hammack, \textsc{Department of Mathematics and Applied Mathematics,\\ 
Virginia Commonwealth University, 
P.O. Box 2014, Richmond, VA  23284 USA}\\
\par\nopagebreak
}}

\maketitle

\begin{abstract}
\noindent 
For strongly connected, pure $n$-dimensional regular CW-complexes, we show that
 {\it evenness} (each $(n{-}1)$-cell is contained in an even number of $n$-cells)
is equivalent to generalizations of both cycle decomposition and traversability.
\end{abstract}


\section{Introduction}
Euler's theorem for connected multigraphs asserts the equivalence of even-degree vertices, decomposition into edge-disjoint cycles, and  existence of a closed trail using all edges.
This basic theorem \cite[p 64]{harary} arose from a problem in recreational mathematics \cite[pp 1--2 ]{harary}, \cite{blw-1976}, and has applications to routing \cite{ej73}, \cite[pp 231--232]{lp86}.  

Recently, Glock, Joos, K\"{u}hn, and Osthus \cite{gjko2020} obtained a higher-dimensional, {\it combinatorial} version of Euler's theorem that applies to regular hypergraphs of a special type, including complete hypergraphs, subject to a finite set of exceptions.  

Here we follow a path pioneered by Gr\"unbaum \cite{branko};
we give higher-dimensional {\it topological} analogues of cycle decomposition and closed Eulerian trail.  This paper is a sequel to \cite{hk-math-mag} (the 2-dimensional case),
which introduced our concepts of {\it circlet} and {\it Euler cover} (full definitions below). 
To extend the argument in \cite{hk-math-mag} from $2$ dimensions to $n$, we use pseudomanifolds  and  
the proof 
requires a key property  (sphericity of codimension-2 intervals in the face poset \cite{bjo},\cite[p 244]{massey-alg-top}) of regular CW-complexes. 

Regular CW-complexes are the spaces obtained by successively attaching closed balls by {\it homeomorphisms} from their boundary (full definitions below). 
The setting of regular CW-complexes is natural since Euler's theorem applies to (and was originally  developed for \cite{euler}) {\it multigraphs} \cite[p 10]{harary}, where multiple 1-cells can be attached to the same pair of vertices but the two boundary vertices of an edge must be distinct. Detailed treatments of regular CW-complexes, including suitability as a framework for pseudomanifolds, are in  \cite[Chap. IX]{massey-alg-top} and \cite{bjo}.

In Section~\ref{Section:Def}, we review regular CW-complexes, describe the concepts, and state our main result
(Theorem \ref{Theorem:Euler}).  Section 3 derives
some properties of our two gadgets, Section 4 proves Theorem~\ref{Theorem:Euler}, and Section 5 surveys the vast literature.  The last section discusses earlier simplicial versions of our gadgets, as well as applications.

\section{Definitions and main result}
\label{Section:Def}
A $k$-{\bf cell} is a closed ball in $\mathbb{R}^k$. A {\bf 0-complex} is a finite set of points.  
In this paper, for $n \geq 1$, we write {\bf n-complex} to denote a finite regular CW-complex \cite[p 94]{massey} 
({\bf regular} means that the $k$-cells, for $1 \leq k \leq n$, are attached by homeomorphisms on their boundaries).
In an $n$-complex $K$, we call the $n$-cells {\bf facets}, the $(n{-}1)$-cells {\bf sides}, and the $(n{-}2)$-cells {\bf corners}. The {\bf degree} $\deg_K(s)$
of a side~$s$ is the number of facets that contain it. 
An $n$-complex  is {\bf even} if each side has even degree and is {\bf pure} if each cell is contained in a facet.   Write $K^{(r)}$ for the set of $r$-cells; the {\bf r-skeleton} $K^r$ denotes the union of all cells
of dimension at most $r$.  The {\bf boundary $\partial(c)$} of an $r$-cell $c$ is an $(r{-}1)$-sphere and $c \setminus \partial(c)$ is the {\bf interior} of $c$.  The interior of a 0-cell is the point itself.  A complex is the disjoint union of the interiors of its cells.  The cells are also called {\bf faces} and the faces define a poset under inclusion.  The topology of a regular CW-complex is determined by its face poset \cite{bjo}.

If $S$ is any non-empty set of $n$-cells of $K$, let {\bf K(S)}, the subcomplex {\bf induced} by $S$, be the intersection of all $n$-subcomplexes $L$ for which $S \subseteq L \subseteq K$. Induced complexes are pure.     
The {\bf dual} {\bf K}$^{{\bf *}}$ is the {\it multigraph} whose vertices are the facets of $K$; an edge joins two vertices once for each side contained in the corresponding two facets. We say $K$
is {\bf strongly connected} (s-connected) if  $K^*$ is connected; $K$ is a {\bf pseudomanifold} if it is pure, s-connected, and each side has degree 2 (\cite{bjo, massey-alg-top, spanier}).  
A pseudomanifold is a {\bf manifold} if each point has a neighborhood homeomorphic~to~$\mathbb{R}^n$.

A continuous function $\varphi: K\to L$ between $n$-complexes 
is a {\bf cellular map} if $\varphi(K^r) \subseteq L^r$ 
for $0 \leq r \leq n$ \cite[p 232]{massey-alg-top}. 

A union of subcomplexes is  facet-disjoint if each facet is in a unique subcomplex.
Write {\bf K} $=$ {\bf J} $\sqcup \,${\bf L}  $\,$ {\bf if K is the facet-disjoint union of J and L.}

One can generalize the concept of cycle by noticing that cycles are
1-manifolds (more specifically, 1-spheres). In higher dimensions, one  
can try to decompose regular CW-complexes into manifolds or spheres. See \cite{hk-bhms, hk-ADAM, hk-geomb, hk-Genus-monthly, pk-tubes-and-cubes} and Section~5.  

A cycle is a {\bf minimal even} graph: (i) its vertices have positive even degree and (ii) there is no proper subgraph whose vertices all have positive even degree. We will call the $n$-dimensional analogue an ``$n$-circlet.''

An $n$-{\bf circlet}  \cite{hk-math-mag} is a pure even $n$-complex which is not the facet-disjoint union of two pure
even $n$-complexes (equivalently, its set of facets cannot be partitioned into
two non-empty sets, each inducing an even $n$-complex); cf. Fogelsanger \cite{fogel}.
\begin{lemma}
\label{lm:even}
For $K$ an even pure $n$-complex,  $K$ is a circlet $\iff$ $K$ is minimal even. 
\end{lemma}
\begin{proof} We show the contrapositives.
If $K$ is a facet-disjoint union of two pure even $n$-complexes 
$K = J \,\sqcup \,L$,
then $J$ is a proper subcomplex of $K$ and $J$ is even.  Conversely, if $J$ is a proper even subcomplex of $K$, then put
$L:= K \setminus J := K( K^{(n)} \setminus J^{(n)})$ so that
$K = J \, \sqcup \,L$. 
To see that $L$ is even, consider side $s \in L^{(n-1)}$.
 If $s \notin J^{(n-1)}$, then $\deg_L(s) = \deg_K(s)$.  If $s \in J^{(n-1)}$, then 
$\deg_K(s) = \deg_J(s) + \deg_L(s)$, where $\deg_K(s) > \deg_J(s)$ and both are positive even, so $\deg_L(s)$ is positive even.  \end{proof}

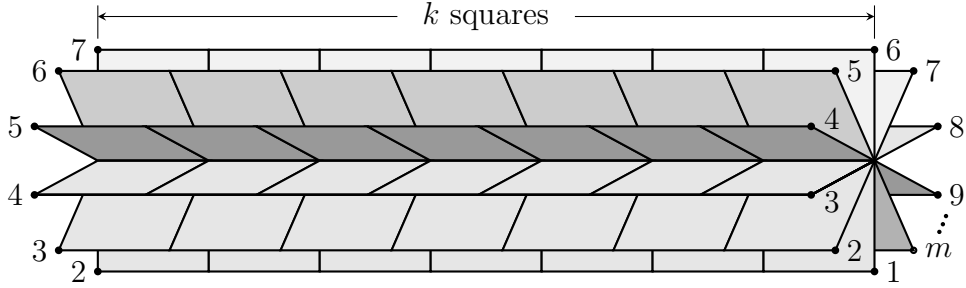
\begin{figure}[t]
\centering
\begin{tikzpicture}[scale=1.475,style=thick]
\def\vr{2.5pt}
\def\sh{7pt}
\def\fac{0.6}

\definecolor{gray0}{rgb}{0.95,0.95,0.95}
\definecolor{gray1}{rgb}{0.9,0.9,0.9}
\definecolor{gray2}{rgb}{0.8,0.8,0.8}
\definecolor{gray3}{rgb}{0.7,0.7,0.7}
\definecolor{gray4}{rgb}{0.6,0.6,0.6}

\def\alpha{18}
\foreach \x in {0,1,2,3,4,5,6} \draw  [fill=gray1, fill opacity=0.9] (\x,0)--({\x+1},0)--({\x+1+\fac*cos(\alpha)},{sin(\alpha)})--({\x+\fac*cos(\alpha)},{sin(\alpha)})--cycle;
\def\alpha{-18}
\foreach \x in {0,1,2,3,4,5,6} \draw [fill=gray4, fill opacity=0.9] (\x,0)--({\x+1},0)--({\x+1+\fac*cos(\alpha)},{sin(\alpha)})--({\x+\fac*cos(\alpha)},{sin(\alpha)})--cycle;

\def\alpha{54}
\foreach \x in {0,1,2,3,4,5,6} \draw  [fill=gray0, fill opacity=0.9] (\x,0)--({\x+1},0)--({\x+1+\fac*cos(\alpha)},{sin(\alpha)})--({\x+\fac*cos(\alpha)},{sin(\alpha)})--cycle;
\def\alpha{-54}
\foreach \x in {0,1,2,3,4,5,6} \draw [fill=gray3, fill opacity=0.9] (\x,0)--({\x+1},0)--({\x+1+\fac*cos(\alpha)},{sin(\alpha)})--({\x+\fac*cos(\alpha)},{sin(\alpha)})--cycle;

\def\alpha{90}
\foreach \x in {0,1,2,3,4,5,6} \draw  [fill=gray0, fill opacity=0.9] (\x,0)--({\x+1},0)--({\x+1+\fac*cos(\alpha)},{sin(\alpha)})--({\x+\fac*cos(\alpha)},{sin(\alpha)})--cycle;
\def\alpha{270}
\foreach \x in {0,1,2,3,4,5,6} \draw [fill=gray1, fill opacity=0.9] (\x,0)--({\x+1},0)--({\x+1+\fac*cos(\alpha)},{sin(\alpha)})--({\x+\fac*cos(\alpha)},{sin(\alpha)})--cycle;
\def\alpha{126}
\foreach \x in {0,1,2,3,4,5,6} \draw  [fill=gray2, fill opacity=0.9] (\x,0)--({\x+1},0)--({\x+1+\fac*cos(\alpha)},{sin(\alpha)})--({\x+\fac*cos(\alpha)},{sin(\alpha)})--cycle;
\def\alpha{234}
\foreach \x in {0,1,2,3,4,5,6} \draw [fill=gray1, fill opacity=0.9] (\x,0)--({\x+1},0)--({\x+1+\fac*cos(\alpha)},{sin(\alpha)})--({\x+\fac*cos(\alpha)},{sin(\alpha)})--cycle;
\def\alpha{162}
\foreach \x in {0,1,2,3,4,5,6} \draw [fill=gray4, fill opacity=0.9] (\x,0)--({\x+1},0)--({\x+1+\fac*cos(\alpha)},{sin(\alpha)})--({\x+\fac*cos(\alpha)},{sin(\alpha)})--cycle;
\def\alpha{198}
\foreach \x in {0,1,2,3,4,5,6} \draw  [fill=gray1, fill opacity=0.9] (\x,0)--({\x+1},0)--({\x+1+\fac*cos(\alpha)},{sin(\alpha)})--({\x+\fac*cos(\alpha)},{sin(\alpha)})--cycle;
\draw (0,0)--(7,0);

\foreach \x in {90,126,162,198,234, 270} \draw ({\fac*cos(\alpha)},{sin(\alpha)})--+(7,0)--(7,0);

\foreach \x in {1,2,3,4,5,6,7,8,9} \draw ({7+\fac*cos(-(\x-9)*36-18)},{sin(-(\x-9)*36-18)}) [fill=black] circle (0.025);
\foreach \x in {1,2,5,6,7,8,9} \draw ({7+\fac*cos(-(\x-9)*36-18)},{sin(-(\x-9)*36-18)}) node [right] {$\x$};
\draw ({7+\fac*cos(158},{sin(158)}) node [right] {$4$};
\draw ({7+\fac*cos(202},{sin(202)}) node [right] {$3$};
\draw ({7+\fac*cos(-54)},{sin(-54)}) circle (0.024)  node [right] {$m$};

\foreach \x in {2,3,4,5,6,7} \draw ({\fac*cos(-(\x-10)*36-18)},{sin(-(\x-10)*36-18)}) [fill=black] circle (0.025) node [left] {$\x$};

\foreach \x in {-30,-35,-40} \draw ({7.125+\fac*cos(\x)},{sin(\x)}) circle (0.01);

\draw [->,>=stealth, thin] (2.7,1.3)--+(-2.7,0);
\draw (3.5,1.3) node {$k$ squares};
\draw [->,>=stealth, thin] (4.3,1.3)--+(2.7,0);
\foreach \x in {0,7} \draw [thin] (\x,1.1)--+(0,0.3);
\end{tikzpicture}
\caption{A family of circlets $C(k,m)$. Let $k  \geq 3$ and let $m \geq 4$ be even.  To construct the 2-circlet $C(k,m)$,
begin with the 2-complex illustrated here, whose 1-skeleton is the graph Cartesian product
$K_{1,m}\Box P_{k+1}$. 
Identify the left $K_{1,m}$ with the right $K_{1,m}$ by a twist of $2\pi/m$, as indicated.
The result is a {\it finned} 2-complex with a $C_{km}$ boundary. Cap this boundary by a polygon with $km$ sides.
The result is a 2-circlet, as the $km$-gon shares an edge of degree 2 with any square, so any even
$n$-subcompex that contains the $km$-gon necessarily contains all the squares too.}
\label{Fig:Circlet}
\end{figure}

Every $n$-pseudomanifold is an $n$-circlet. For $n=1$ every circlet is a pseudomanifold (in fact, a sphere), but for $n>1$ that is no longer the case, see Figure~\ref{Fig:Circlet}. 

A continuous map  $\varphi:M\to K$ is an {\bf Euler cover} if for $n \geq 1$, $K$ is a pure even $n$-complex, $M$ is an $n$-pseudomanifold, $M^{n{-}2} =K^{n{-}2}$,  $\varphi|_{M^{n{-}2}} = \mbox{id}$,
and $\varphi$ induces a bijection from $M^{(n)}$ to $K^{(n)}$ mapping facets homeomorphically. Note $\varphi$ is cellular.


The case $n=1$ deserves special mention. As a 1-pseudomanifold is just a cycle, $M$ is a cycle, and  $M^{n{-}2} = K^{n{-}2}=\emptyset$, so
the condition $\varphi|_{M^{n{-}2}} = \mbox{id}$ holds vacuously.
Figure~\ref{Fig:EC1} shows a 1-dimensional Euler cover, and illustrates that this coincides with the 
notion of an Euler circuit. Figure~\ref{Fig:K6Cover} gives an example of a 2-dimensional Euler cover.
(In \cite{hk-math-mag}, we addressed the case $n=2$ and required $M$ to be a {\it manifold}.)

We can now formulate our result.



\begin{theorem}
\label{Theorem:Euler}  
Let $K$ be a pure, s-connected $n$-complex. The following are equivalent:\begin{enumerate}[itemsep=-5pt, align=left,topsep=2pt, labelwidth=1.5em,leftmargin=1.5cm]\item [{\em (i)}] $K$ is even,
\item [{\em (ii)}] $K$ is a facet-disjoint union of circlets,
\item [{\em (iii)}] There exists an $n$-pseudomanifold $M$ and an Euler cover $\varphi: M \to K$.
\end{enumerate}\end{theorem}

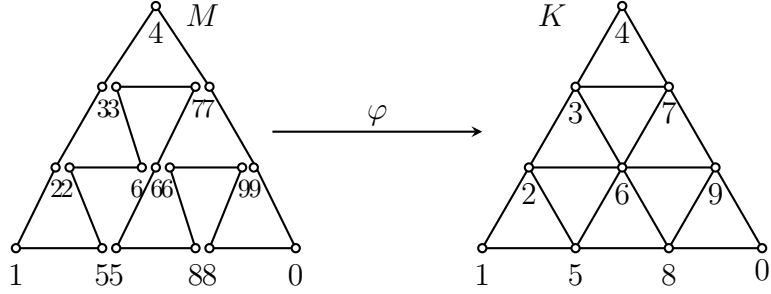
\begin{figure}[h]
\centering
\begin{tikzpicture}[scale=0.62,style=thick]
\def\vr{2.5pt}
\def\sh{7pt}

\draw (0,0)--+(60:6)--+(6,0)--cycle;
\draw (60:2)--+(4,0) (60:4)--+(2,0);
\draw (2,0)--+(60:4) (4,0)--+(60:2);
\draw (2,0)--+(120:2) (4,0)--+(120:4);
\foreach \x  in {1,2,3,4} \draw (60:{2*(\x-1)}) [fill=white] circle (\vr)  +(0,-0.1) node [below] {\raisebox{3pt}{$\x$}};
\foreach \x  in {5,6,7} \draw ({2+2*(\x-5)*cos(60)}, {2*(\x-5)*sin(60)}) [fill=white] circle (\vr) +(0,-0.1) node [below] {\raisebox{3pt}{$\x$}};
\foreach \x  in {8,9} \draw  ({4+2*(\x-8)*cos(60)}, {2*(\x-8)*sin(60)}) [fill=white] circle (\vr)  +(0,-0.1) node [below] {\raisebox{3pt}{$\x$}};
\draw (6,0)  [fill=white] circle (\vr)  node [below] {\raisebox{3pt}{$0$}};
\draw (1.5,5) node {$K$};

\draw [->,>=stealth] (-4.5,2.5)--(0,2.5);
\draw (-2.25,2.85) node {$\varphi$};

\draw (-7,3*1.732)--(-8.15,2*1.732)--(-9.15,1*1.732)--(-10,0)--(-8.15,0)--(-8.85,1*1.732)--(-7.3,1*1.732)--(-7.85,2*1.732)--(-6.15,2*1.732)--(-7.85,0)--(-6.15,0)--(-6.7,1*1.732)--(-5.15,1*1.732)--(-5.85,0)--(-4,0)--(-4.85,1.732)--(-5.85,2*1.732)--(-7,3*1.732);
\foreach \x  in {0,1.85,2.15,3.85,4.15,6} \draw (-10,0) +(\x,0) [fill=white] circle (\vr);
\foreach \x  in {-0.15, 0.15, 1.7, 2,2.3, 3.85,4.1} \draw (-9,0) +(\x,1.732) [fill=white] circle (\vr);
\foreach \x  in {-0.15, 0.15,  1.85, 2.15} \draw (-8,0) +(\x,2*1.732) [fill=white] circle (\vr);
\draw (-7,3*1.732) [fill=white] circle (\vr);
\draw (-6,5) node {$M$};
\draw (-10,-0.1) node [below] {1};
\draw (-8.15,-0.1) node [below] {5};
\draw (-7.85,-0.1) node [below] {5};
\draw (-6.15,-0.1) node [below] {8};
\draw (-5.85,-0.1) node [below] {8};
\draw (-4,-0.1) node [below] {0};
\draw (-7,3*1.732-0.1) node [below] {4};
\draw (-6.975,1*1.732-0.03) node [below] {\footnotesize 6};
\draw (-7.4,1*1.732-0.03) node [below] {\footnotesize 6};
\draw (-6.75,1*1.732-0.03) node [below] {\footnotesize 6};
\draw (-9.125,1*1.732-0.03) node [below] {\footnotesize  2};
\draw (-8.9,1*1.732-0.03) node [below] {\footnotesize  2};
\draw (-5.1,1*1.732-0.03) node [below] {\footnotesize  9};
\draw (-4.875,1*1.732-0.03) node [below] {\footnotesize  9};
\draw (-8.125,2*1.732-0.03) node [below] {\footnotesize  3};
\draw (-7.9,2*1.732-0.03) node [below] {\footnotesize  3};
\draw (-6.1,2*1.732-0.03) node [below] {\footnotesize  7};
\draw (-5.875,2*1.732-0.03) node [below] {\footnotesize  7};
\end{tikzpicture}
\caption{An Euler cover of an even 1-complex $K$.  A vertex labeled $x$ in $M$ maps to $x$ in  $K$.
An edge labeled $xy$ in $M$ maps to $xy$ in $K$.
}
\label{Fig:EC1}
\end{figure}
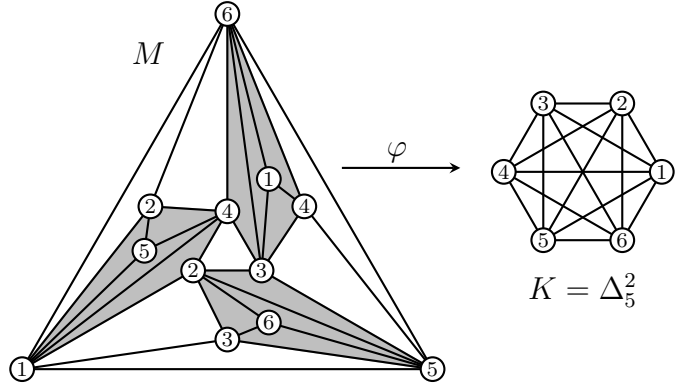
\begin{figure}[h]
\centering
\begin{tikzpicture}[scale=0.525,style=thick]
\def\h{6}
\def\hh{1}
\def\hhh{2.25}
\def\vr{0.3}
\def\sh{4.5}
\def\tr{2.1} 
\draw [color=lightgray, fill=lightgray, fill opacity=0.5] (90:\h)--(30:\hhh)--(330:\hh)--(90:\hh)--(90:\h);
\draw [color=lightgray, fill=lightgray, fill opacity=0.5] (90+120:\h)--(30+120:\hhh)--(330+120:\hh)--(90+120:\hh)--(90+120:\h);
\draw [color=lightgray, fill=lightgray, fill opacity=0.5] (90-120:\h)--(30-120:\hhh)--(330-120:\hh)--(90-120:\hh)--(90-120:\h);
\foreach \x  in {90,210,330} \draw (\x:\hh)--(\x+120:\hh);
\foreach \x  in {90,210,330} \draw (\x:\h)--(\x+120:\h);
\foreach \x  in {-90,30,150} \draw (\x:\hhh)--(\x+60:\h);
\foreach \x  in {-90,30,150} \draw (\x:\hhh)--(\x-60:\h);
\foreach \x  in {-90,30,150} \draw (\x:\hhh)--(\x-60:\hh);
\foreach \x  in {90,210,330} \draw (\x:\hh)--(\x+120:\h);
\foreach \x  in {90,210,330} \draw (\x:\hh)--(\x:\h);
\foreach \x  in {180,60,300} \draw (\x-30:\hhh)--(\x:\tr)--(\x+30:\h);
\foreach \x  in {90,210,330} \draw (\x:\hh)--(\x+90:\tr);
\foreach \x  in {90,210,330} \draw (\x:\h) [fill=white] circle (\vr);
\foreach \x  in {90,210,330} \draw (\x:\hh) [fill=white] circle (\vr);
\foreach \x  in {-90,30,150} \draw (\x:\hhh) [fill=white] circle (\vr);
\foreach \x  in {180,60,300} \draw (\x:\tr) [fill=white] circle (\vr);
\draw (90:\h) node {\scriptsize 6};
\draw (210:\h) node {\scriptsize 1};
\draw (330:\h) node {\scriptsize 5};
\draw (90:\hh) node {\scriptsize 4};
\draw (210:\hh) node {\scriptsize 2};
\draw (330:\hh) node {\scriptsize 3};
\draw (-90:\hhh) node {\scriptsize 3};
\draw (30:\hhh) node {\scriptsize 4};
\draw (150:\hhh) node {\scriptsize 2};
\draw (180:\tr) node {\scriptsize 5};
\draw (60:\tr) node {\scriptsize 1};
\draw (300:\tr) node {\scriptsize 6};
\draw [->,>=stealth] (2.9,2.1)--+(3,0);
\foreach \x in {1,2,3,4,5,6} \draw (9,2) +(\x*60-60:2)--+(\x*60:2);
\foreach \x in {1,2,3,4,5,6} \draw (9,2) +(\x*60-60:2)--+(\x*60+60:2);
\foreach \x in {2,4,6} \draw (9,2) +(\x*60:2)--+(\x*60+180:2);
\foreach \x in {1,2,3,4,5,6} \draw (9,2) +(\x*60-60:2) [fill=white] circle (\vr) node {\scriptsize $\x$};
\draw (-2,5) node {$M$};
\draw (9,-1) node {$K=\Delta_5^{2}$};
\draw (4.3,2.5) node {$\varphi$};
\end{tikzpicture}
\caption{An Euler cover $\varphi\colon M\to \Delta_5^{2}$ of the 2-skeleton of the
5-dimensional simplex $\Delta_5$, which has six vertices $1,2,3,4,5,6$, fifteen edges and twenty triangular faces.
The pseudomanifold $M$ has six vertices, thirty edges and twenty triangular faces (including the unbounded
face in this planar drawing).
Vertices with the same labels are identified, so $M$ is a sphere with six pinchpoints.
Note that the three shaded areas of $M$ map to three pairwise face-disjoint tetrahedron boundaries in $K$. 
The white triangles (including the unbounded region) map to the boundary of an octahedron in $K$.
This decomposes $K$ as a face-disjoint union of four spheres (circlets).
} 
\label{Fig:K6Cover}
\end{figure}

\section{Some properties of circlets and pseudomanifolds}
We begin with some lemmas that will be needed in the proofs.
The following is stated in \cite[p 244]{massey-alg-top} (see 
\cite{bjo} for the stronger property ``interval sphericity'' mentioned above).  
\begin{lemma}
In an $n$-complex, for $2 \leq r \leq n$, if an $r$-cell $\gamma$ contains an $(r{-}2)$-cell $\alpha$, then there are exactly two $(r{-}1)$-cells $\beta_1$, $\beta_2$ with $\alpha \subset \beta_j \subset \gamma$, $j = 1,2$.
\label{lm:massey-244}
\end{lemma}
\begin{proof}
The interiors of the cells partition the complex, so $\alpha \subset \partial{\gamma}$.  But $\gamma$ is an $r$-ball so its boundary is an $(r{-}1)$-sphere, which is a manifold and hence a pseudomanifold.  Thus, $\alpha$ is in the boundary of exactly two $(r{-}1)$-cells, $\beta_1$ and $\beta_2$, which are both in $\partial \gamma$ and so contained in $\gamma$.
\end{proof}
Klee \cite{klee} proved that the dual graph $M^*$ of a simplicial $n$-manifold $M$ is $(n{+}1)$-connected, but this depends on the simplicial structure; the dual of a regular CW-complex can have vertices of degree 2.  In a connected multigraph, an edge is a {\bf bridge} if its removal disconnects the graph. 
The following analogue of Klee's result uses Lemma \ref{lm:massey-244} and is applied to prove strong connectedness in Proposition 2 below. 
\begin{lemma}
Let $M$ be any regular CW-complex that is an $n$-pseudomanifold.  Then the dual multigraph $M^*$ is bridgeless.
\label{th:bridgeless}
\end{lemma}
\begin{proof}
We show that any edge $e$ of $M^*$ lies on a cycle. Suppose $e$ is determined by side $s$.
Since $M$ is a pseudomanifold, $s$ lies in exactly two facets $f_1, f_2 \in M^{(n)}$. Choose a corner $c$ of $M$ that is contained in $s$ and let $\cF(c)$ be the set of all facets of $M$ that contain~$c$. 
Observe that 
$f_1, f_2 \in \cF(c)$.  
Now define a multigraph $M^*(c)$ with vertex set $\cF(c)$, where facets $f$ and $f'$ are adjacent in $M^*(c)$ once for every side $s$ such that $c \subset s$ and $s \subseteq \partial f \cap \partial f'$. 
By Lemma \ref{lm:massey-244}, each vertex has degree 2. Hence, $M^*(c)$ is a vertex-disjoint union of cycles.  But $e$ is an edge of $M^*(c) \subseteq M^*$. 
\end{proof}

An interesting aspect of circlets (though not needed in the development below) is that they form a matroid. Indeed, a {\bf circuit matroid} \cite[p 40]{harary} is a given finite set $X$ and a family $\mathcal{C}$ of non-empty subsets of $X$ such that\\ 
(a) no proper subset of an element in $\mathcal{C}$ is in $\mathcal{C}$, and\\ 
(b) if $C_1, C_2 \in \mathcal{C}$, $C_1 \neq C_2$, and $x \in C_1 \cap C_2$, then there is an element $C \in \mathcal{C}$ such that $C \subseteq C_1 \cup C_2 \setminus \{x\}$.  
\begin{lemma}
The circlets of an $n$-complex $K$ form a circuit matroid on $K^{(n)}$.
\label{lm:circ-matroid}
\end{lemma}
\begin{proof}
Let $\mathcal{C}$ be the collection of all subsets of $K^{(n)}$ which induce circlets.  By our Lemma \ref{lm:even},
condition (a) holds for $\mathcal{C}$. To get (b), note that for $C_1 \neq C_2 \in \cC$ and $x$ any facet in $C_1 \cap C_2$, we have $C_1 \,\Delta \, C_2 \subseteq C_1 \cup C_2 \setminus \{x\}$, where
\[
C_1 \,\Delta \, C_2 := (C_1 \setminus C_2) \cup (C_2 \setminus C_1) =
(C_1 \cup C_2) \setminus (C_1 \cap C_2) \subseteq (C_1 \cup C_2) \setminus \{x\}.
\]
If $s \in K(C_1 \,\Delta \, C_2)^{(n-1)}$, then $t := \deg(s,K(C_1 \,\Delta \, C_2)) > 0$. 
If $s \in K(C_1 \setminus C_2)^{(n-1)}$, then $t = \deg(s,C_1) - \deg (s,C_2)$ which is the difference of even numbers. Hence, $K(C_1 \,\Delta \, C_2)$ is an even complex and contains
some member of $\mathcal{C}$ by Lemma \ref{lm:even}.
\end{proof}





\section{
Machinery and proof of main theorem}
\label{Section:chec}

Here we introduce notation, construct Euler covers,  and prove Theorem \ref{Theorem:Euler}.

Let $K$ be a pure, even $n$-complex. 
For each $s$ in $K^{(n-1)}$, let $\mathcal{F}_s  \subseteq K^{(n)}$ be the set of facets of $K$ that contain side $s$.  Let $\mathcal{P}_s := (P_1, \ldots, P_{r_s})$ be a fixed ordered partition of $\mathcal{F}_s$ into 
$$r_s := | \mathcal{F}_s|/2 = \deg_K(s)/2$$ {\it unordered pairs} $P_j$ of distinct facets.  
When necessary, we show the complex from which pairs or facets are taken by a superscript (e.g., if $L$ is a subcomplex of $K$ and $s$ is a side of $L$, one has $\mathcal{F}^L_s \subseteq \mathcal{F}^K_s $).
Given $\sigma\in\mathcal{F}_s$, let $g_s(\sigma)$ 
satisfy $\{\sigma, g_s(\sigma)\} = P_j$, for some $j$, $1 \leq j \leq r_s$.
Observe that $g_s$ is an involution on $\mathcal{F}_s$ for each side $s$ of $K$. 

Write $g$ for the 
collection of involutions $\{g_s: s \in K^{(n-1)}\}$ arising from a set of partitions $\{\mathcal{P}_s:s \in K^{(n-1)}\}$.
(Informally $g$ is referred to as a ``gluing,'' because it will be used to attach a facet $\sigma$ to facet $g_s(\sigma)$
in a manner described in the proof of Proposition~\ref{Prop:CircletCover}.)
For each side $s \in K^{(n-1)}$ and facet $\sigma \in \mathcal{F}_s$, let
$\nu(s,\sigma)$ be the unique index $1 \leq \nu(s,\sigma)\leq r_s$
for which 
\begin{equation}
\label{Prop:CircletCover}
\mathcal{P}_s = (P_1, \ldots , P_{r_s}) \; \implies P_{\nu(s,\sigma)}= \{\sigma, g_s(\sigma)\}.
\end{equation} 
\begin{proposition}
Every circlet $K$ has an Euler cover. 
\label{pr:crclt-cover}
\end{proposition}
\begin{proof}
We first build an $n$-complex $M$  from $K$ and any ``gluing'' $g$ as follows:
\smallskip

\n 
 (1) 
 Define $M^{n-2} := K^{n-2}$ (there is no change in the $(n{-}2)$-skeleta!);\\
 (2) 
Put $M^{(n-1)} := \{s \times \{j\}: s \in K^{(n-1)}, 1 \leq j \leq r_s \}$ where each $s \times \{j\}$ is a copy of $s$ attached to $M^{(n-2)}$
in the same way as  $s$ is attached to $K^{(n-2)}$. Thus, $M^{n-1}$ is obtained from $K^{n-1}$ by removing the interior of each side $s$ and replacing it with $r_s$ identical copies (like parallel edges) all with boundaries attached identically to the $(n{-}2)$-skeleton of $M$ as $s$ is attached to the corresponding skeleton of $K$.\\
\noindent (3)  
{\it This is where the gluing $g$ is used.}
Let $M^{(n)} := K^{(n)}$ but we will attach these $n$-cells differently and denote them  $ \{\hat{\sigma}: \sigma \in K^{(n)}\}$.
For each facet $(\sigma,a)$ of $K$, with attaching map $a: \partial{\sigma} \to K^{n-1}$, 
$a\big(\partial{\sigma}\big)$ is an $(n{-}1)$-sphere in $K^{n-1}$ and 
\[
a(\partial{\sigma})= (a(\partial{\sigma}))^{n-2} \cup \bigcup_{s \in a(\partial{\sigma})^{n-1}} {\rm int}(s)
\]
is the disjoint union of its $(n{-}2)$-skeleton and the interiors of its $(n{-}1)$-cells.  Using the same partition of $\partial{\hat{\sigma}}$, we define 
 $$\hat{a}: \partial{\hat{\sigma}} \to M^{n-1}$$
 to be the identity on the $(n{-}2)$-skeleton but on the interior of the $(n{-}1)$-cell in $\partial{\hat{\sigma}}$ corresponding to $s$, instead of attaching to $s$, one has 
$$\hat{a}|_{a^{-1}(s)} =  a|_{a^{-1}(s)} \times \{\nu(s,\sigma)\}.$$ 
The map $\hat{a}$ is a homeomorphism as it is the identity on the $(n{-}2)$-skeleton, a homeomorphism on the disjoint interiors of the sides, and all $s \times \{j\}$ have disjoint interiors.
\\  
\vspace{-.4cm}

{\it Define a continuous map} $\varphi: M \to K$ which extends the identity on the $(n{-}2)$-skeleta by mapping $s \times \{j\}$ to $s$ by first-coordinate projection $(x,j) \mapsto x$ for $x \in s$ (for each $s$, $1 \leq j \leq r_s$), and mapping $\hat{\sigma}$ by the identity on the interior of the facet.  

We now check that $M$ is a pseudomanifold.
Indeed, each side is of degree 2. Since $K$ is pure, so is $M$. 
To see that $M$ is s-connected,
let $N$ be a {\bf strong component}  of $M$; that is, let $N$ be the subcomplex of $M$ corresponding to a connected component of $M^*$.   
Let $L:= K\big(\varphi(\{\sigma : \sigma \in N^{(n)} \})\big)$ denote the subcomplex of $K$ induced by the $\varphi$-images of $N^{(n)}$.  
Any side $s$ in $L$ is the image of $r(N,s) := |\{j \in [r_s]: s \times {j} \in N\}|$ sides in $N$ so $s$ is incident to the $\varphi$-images of the $2 r(N,s)$ facets of $N$ incident to the sides in $N$ mapping to $s$.  Thus, $L$ is even and by Lemma \ref{lm:even}, $L=K$, so $N = M$; hence, $M$ is a pseudomanifold, and $\varphi: M \to K$ is an Euler cover. 
\end{proof}

Put $\widehat{K_g} := M$ and
let $\varphi_g: \widehat{K_g} \to K$ be as in the proof of Proposition \ref{pr:crclt-cover}.  The map $\kappa: g \mapsto \varphi_g$ is an {\it injection} from the set of gluings of the circlet $K$ to the set of its Euler cover maps.   
But the pseudomanifolds themselves can be homeomorphic. (There are very many gluings but in dimension 2, the feasible range of the Euler characteristic of the resulting surface isn't large \cite{hk-math-mag}).
The correspondence is also {\it surjective}. For any Euler cover $\varphi: V \to K$, the gluing $g_{\varphi}$ of $K$, defined for $s \in K^{(n-1)}$ by $(g_\varphi)_s(\sigma) = \sigma'$ if and only if $\{\sigma,\sigma'\}$ is the image under $\varphi$ of two faces in $V$ meeting in a side $\hat{s}$ for which $\varphi(\hat{s})=s$, satisfies $\varphi = \varphi_{(g_\varphi)}$. 
 
One also may deduce that for any Euler cover $\varphi:V \to K$:
for each $\hat{s} \in V^{(n-1)}$, $\varphi(\hat{s})=s \in K^{(n-1)}$; for each $s \in K^{(n-1)}$, $\varphi^{-1}({\rm int}(s)) = {\rm int}(s) \times \{1,\ldots, r_s\}$, where $r_s = \deg_K(s)/2$; and $\varphi|_{s \times \{j\}}$ is first-coordinate projection for $1 \leq j \leq r_s$. Hence, an Euler cover must map each cell homeomorphically to a cell of the same dimension, but the correspondence of sides in $V$ to sides in $K$ is a many-to-one surjection.

For $\sigma, \tau, \mu, \lambda \in \cF^K_s$ with $\{\sigma,\tau\}$ and $\{\mu,\lambda\}$ paired by a gluing, a {\bf recombination} is a pairing $\{\sigma,\mu\}, \{\tau,\lambda\}$ or $\{\sigma,\lambda\}, \{\tau,\mu\}$.
Recall ``$\sqcup$'' means facet-disjoint union.
\begin{proposition}
Let $K = L \,\sqcup \,C$ be a strongly connected pure $n$-complex, let   
$L$ have an Euler cover, and let $C$ be a circlet. 
Then $K$ has an Euler cover. 
\label{prop:induct}
\end{proposition}
\begin{proof}
Since $K = L \,\sqcup \,C$ is strongly connected, there exists some $s$ in $L^{(n-1)} \cap C^{(n-1)}$.
Let $\psi: N \to L$ be an Euler cover. By Proposition \ref{Prop:CircletCover} and its proof, there is an Euler cover $\xi: Y \to C$ with $N \cap Y = N^{(n-2)} \cap Y^{(n-2)} = L^{(n-2)} \cap C^{(n-2)}$, so 
we may choose a {\it distinct} pair of sides $\hat{s}_i$, $\tilde{s}_j$, with $\hat{s}_i \in \psi^{-1}(s)$ for some $1 \leq i \leq r_s^L$ and  $\tilde{s}_j \in \xi^{-1}(s)$, for $1 \leq j \leq r_s^C$.  As $L$ and $C$ are face-disjoint, ${\cal F}_s^K = {\cal F}_s^L \cup {\cal F}_s^C$ (disjoint union).

We first sketch the remainder of the argument and then fill in the details.
 Apply recombination to $P_i^L$ and $P_j^C$, but keep all other facet-pairs, determined by $\psi$ and $\xi$, resp.,  unchanged.  Let $g$ be the resulting gluing and let $M = \widehat{K}_g$. Then $M$ is a pseudomanifold and there is an Euler cover from $M$ to $K$. The details now follow.

There are two unique facets $\sigma'$, $\tau'$ in $N$ containing $\hat{s}_i$ and two unique facets $\mu'$, $\lambda'$ in $Y$ containing $\tilde{s}_j$.  Let $\sigma, \tau, \mu, \lambda$ be the cells in $L$ and $C$, resp., to which $\sigma', \tau', \mu', \lambda'$ project under $\psi$ and $\xi$, resp.  
Then $P_i^L = \{\sigma,\tau\}$ and $P_j^C = \{\mu,\lambda\}$. 

Let ${\cal P}_s^L * {\cal P}_s^C$ be the  concatenation of the two sequences of unordered facet-pairs and let ${\cal P}_s^K$ be the result of applying a recombination to it
so that $P_i^L$ and $P_j^C$ are deleted while $\{\sigma, \mu\}$ and $\{\tau, \lambda\}$ are added in their place.  

For any side $t$, ${\cal P}_t^K := {\cal P}_t^L *  {\cal P}_t^C$; one of the sequences is empty unless $t$ is in both $L$ and $C$.  
This defines a family of pairwise-partitions $\{ P_s^K: s \in K^{(n-1)}\}$ and so a gluing $g$.  
The regular CW-complex $M := \widehat{K_g}$ is obtained from $N \cup Y$ by 
changing the attaching map of the $n$-cell $\tau'$ so it attaches to $\tilde{s}_j$ instead of $\hat{s}_i$ and of $\mu'$ so it attaches to $\hat{s}_i$ (instead of $\tilde{s}_j$).  
All other attaching maps for cells in  $M$ remain the same.
By construction, $M$ is pure and every side has degree 2.
Further, $M$ is  strongly connected. Indeed, $M^* = (N^* - \sigma' \tau')  \cup(Y^* - \mu' \lambda') + \sigma' \mu' + \tau' \lambda'$.  The first and second parenthesized terms are connected by Lemma \ref{th:bridgeless} so adding the last two edges $\sigma' \mu'$ and $ \tau' \lambda'$, we obtain a connected multigraph. Thus, $M$ is a pseudomanifold.

An Euler cover $\varphi: M \to K$ is now determined by the cellular maps $\psi$ and $\xi$, and the recombination. As $N$ and $Y$ are disjoint in dimensions $n{-}1$ and $n$, while being identical in lower dimensions, the functions $\psi$ and $\xi$ define a cellular map
\begin{equation}
\psi \cup \xi: N  \cup Y \to L \,\sqcup \,C,
\end{equation}
and $\psi \,\cup \,\xi$ is the identity on the $(n{-}2)$-skeleton.  Indeed, the interiors of the $n$ and $(n{-}1)$-cells are pairwise-disjoint, so continuous functions defined on open cells in $N$ and open cells in $Y$ define a continuous function on the union.  
Moreover, $\psi$ and $\xi$ induce bijections 
$\psi_n: N^{(n)} \to L^{(n)}\; \;\mbox{and}\; \;\psi_n: Y^{(n)} \to C^{(n)}$,
so $(\psi \,\cup\,\xi)_n: (N \cup Y)^{(n)} \to K^{(n)}$
is a bijection. 
Let $\varphi: M \to K$  be the continuous map induced by $\psi \,\cup \,\xi$. 
Then $\varphi$ is an Euler cover of $K$.
\end{proof}

We can now prove Theorem~\ref{Theorem:Euler}. Let $K$ be a pure, strongly connected $n$-complex.

\begin{proof}
(i)$\Rightarrow$(ii) (even implies decomposable). 
Let $K$ be even. If it cannot be decomposed, 
$K$ is a circlet.
If $K=K_1 \, \sqcup \, K_2$, where $K_1$, $K_2$ are even, then each of $K_1$ and $K_2$
is either a circlet or can be decomposed as a face-disjoint union of even $n$-complexes. 
 By finiteness, $K=K_1 \sqcup \cdots \sqcup K_n$ is a facet-disjoint union of circlets.

(ii)$\Rightarrow$(iii) (decomposable implies coverable). Let $K$ be
a facet-disjoint union of $r$ circlets.  
If $r=1$, use Proposition~\ref{Prop:CircletCover}. For  $r \geq 2$, form a graph with vertices the set of $r$ circlets, where two facet-disjoint circlets are adjacent if they intersect in a side of $K$.  As $K$ is strongly connected, this graph is connected.  Any nontrivial connected graph has at least two non-cutpoints, corresponding to vertices at diametric distance. Hence, there exists a circlet $C$ in $K$ such that $K = K' \sqcup C$, where $K'$ is both strongly connected and the facet-disjoint union of $r-1$ circlets.  By Proposition~\ref{prop:induct}, truth for $r$ follows from truth for $r-1$, which holds inductively.

(iii)$\Rightarrow$(i) (coverable implies even). Let $\varphi: M \to K$ be an Euler cover. Each side of $M$ is in exactly two facets, so the degree of a side in $K$ is twice the number of sides in $M$ which map to it.
\end{proof}


\section{Historical Context}
\label{Section:Intro}
Euler saw that the following are equivalent for a connected multigraph $G$ \cite[p 64]{harary}:
 
 \smallskip
\begin{enumerate}[itemsep=-4pt, align=left,topsep=2pt, labelwidth=1.5em,leftmargin=1.5cm]
\item [(i)] (evenness) Each vertex of $G$ has positive even degree,
\item [(ii)] (decomposability) $G$ is an edge-disjoint union of cycles,
\item [(iii)] (traversability) $G$ has an Eulerian tour (a closed trail covering each edge).
\end{enumerate}

\smallskip
  
Though the theorem was implicit in
Euler's paper  \cite{euler, euler2} from 1736,
rigorous proof of the equivalence of evenness and decomposability-into-cycles was only first given by 
Veblen \cite{veblen-1912} in 1912, while proof that evenness is equivalent to the existence of a closed Eulerian trail is due to Hierholzer (1873) \cite{hierholzer-1873}; 
see Biggs et al. \cite[pp. 1--11]{blw-1976} for a translation of  \cite{euler} from the Latin.  

We now review the literature.  There are five distinct parts: (1) Our previous work motivated by an attempt to extend the equivalence of (i) and (ii); (2) Earlier work in  similar directions also implicitly aimed at an Eulerian equivalence; (3) The general program \cite{branko} of replacing graphs by complexes; (4) The work on hypergraphs which aimed at an Euler equivalence between (i) and (iii), replacing each by appropriate notions for higher order hypergraphs; (5) Generalizations of (iii) to hypergraphs.

(1) In \cite[p 533]{hk-GC} we proposed to generalize Euler's theorem by decomposing even $2$-dimensional convex cell complexes into facet-disjoint unions of manifolds.  As a step toward this goal, we showed that the even 2-skeleta of the {\it Platonic} complexes can be split into tori and spheres \cite{hk-bhms, pk-tubes-and-cubes} and managed with only spheres (for most cases) in \cite{hk-ADAM}.  Further, we provided facet-disjoint sphere decompositions in \cite{hk-geomb} for even $k$-skeleta of $n$-dimensional Platonic polytopes for all $k$, both constructively (in some special cases) and, in general, asymptotically and existentially (using Keevash \cite{keevash}).  We also partitioned the 2-skeleton of the (odd-dimensional) $d$-hypercube polytope ($d \geq 3$ odd) in \cite{hk-Genus-monthly} into copies of the genus surface of the hypercube graph $Q_d$. 

(2) The idea of finding surfaces within complexes and of  decomposing skeleta of complexes into face-disjoint unions of surfaces had appeared earlier (as we discovered).
In his 1974 thesis \cite{schulz-thesis}, Schulz defined {\bf r-Hamiltonian} $k$-manifolds in the $k$-skeleton of a complex $K$ as subcomplexes of $K$ which are closed manifolds and contain the entire $r$-skeleton $K^r$; see also \cite{bsw76, k-book, ks-1991} and especially Coxeter \cite{cox-skew}. 

Spreer \cite{spreer} observes that 
one may want to write the skeleton of a triangulated manifold as a facet-disjoint union of parts with given properties,
e.g., compact manifolds or compact manifolds with boundary.  
One could use the closed simplexes (or cells), but might prefer a smaller number of parts.  For the
$2$-skeleton of the $d$-cross-polytope, he found \cite{spreer} a set of facet-disjoint $0$-Hamiltonian tori and Klein bottles while we partitioned it into tetrahedral and octahedral boundary spheres \cite{hk-bhms}.

A related issue is {\it Which surfaces are homeomorphic to subcomplexes of polytopes?}
Altshuler  \cite{altshuler71DM} showed that  the sphere and the torus are the only surfaces which are subcomplexes of a {\it cyclic} {\rm 4-}polytope, 
and the torus must $0$-Hamiltonian.  In contrast, a {\it stacked} {\rm 4-}polytope \cite{altshuler71JCT},
has a much larger family of embeddable surfaces .  
Effenberger and K\"{u}hnel \cite{ek-2010} list the  $1$-Hamiltonian surfaces for the Platonic polytopes. 

(3) Sos, Erdos, and Brown \cite{seb73} solved a problem for $2$-complexes, analogous to Turan's:
``{\it What is the maximum number of $2$-simplexes a simplicial-$2$ complex may contain without containing a subcomplex which is a triangulation of the $2$-sphere}?''
If the complex has $n$ vertices, the answer is O$(n^{3/2})$.
The idea of {\it spanning tree of a graph} has been generalized to complexes by  Pippert and Beineke \cite{pippert-beineke}, Dewdney \cite{dewdney-trees}, and Bolker \cite{bolker}, and studied as {\it cellular trees} (e.g., Duval, Klivans, \& Martin \cite{dkm-2016}) and {\it hypertrees} (e.g., 
Kalai \cite{kalai-hyper}, 
Mathew, Newman, Rabinovich, \&  Rajendraprasa \cite{mathew-2018}). 

(4) The generalization of Glock-et-al \cite{gjko2020} 
is based on 
the idea of ``universal cycle'' by Chung, Diaconis, and Graham \cite{cdg-89, cdg-92}.  
A cyclic sequence of vertices in an $r$-uniform hypergraph is a {\bf tight Eulerian tour} \cite{gjko2020} if each (moving) window of length $r$ is an $r$-set and each hyper-edge appears exactly once.  Thus, {\it ``Eulerian tour'' retains the property of being a closed sequence of vertices} in their hypergraph model.

It follows from existence of a tight Eulerian tour that $r$ divides the vertex degrees; conversely, \cite{cdg-89, cdg-92} conjectured that for large $n$, divisibility implies existence. 
This was proved in \cite{gjko2020} using a method of \cite{gklo-2020} which extends \cite{keevash} (existence of designs) from complete $r$-uniform hypergraphs to a larger set of $r$-uniform hypergraphs (with sufficiently many vertices) satisfying  certain conditions.
A {\bf tight r-uniform m-cycle} is a cyclic sequence of $m$ vertices such that each $r$-term consecutive subsequence appears exactly once.  For $m = 2r$, Glock-et-al. \cite{gjko2020} show that, again under certain technical conditions in addition to divisibility, a decomposition into such tight cycles always exists for $m$ sufficiently large.

(5) The concepts of {\it Eulerian} and {\it quasi-Eulerian hypergraph} are considered in Wagner's 2019 thesis \cite{wagner-2019} and in 
Bahamian \& Sajna  \cite{bs-2017} and Sajna \& Steimle \cite{ss-2018}.
Bermond et al. \cite{bermond} introduced {\it Hamiltonian hypergraphs} according to 
\cite{kat-kier}. 

For an $n$-complex $K$, let $\cH(K)$ be the non-uniform hypergraph with $V(\cH(K)) = K^{(n)}$, where the hyperedges are the $\cF_s$, $s \in K^{(n-1)}$.   
If $K$ is Eulerian, we wonder if the hypergraph $\cH(K)$ is Eulerian or Hamiltonian. 
Perhaps this is related to ``Eulerian families'' in hypergraphs  \cite{bs-2017, ss-2018, wagner-2019}.

\section{Discussion}

The notions of evenness and of circlet have also appeared previously in the literature.  Fogelsanger \cite{fogel} defined a minimal simplicial homology cycle over an Abelian group in his 1988 thesis. Our notion of circlet corresponds to the special case of $\mathbb{Z}_2$-coefficients but is also more general since it applies to regular CW-complexes. See also  Cruickshank, Jackson \& Tanigawa \cite[p 3]{cjt} who call simplicial circlets ``simplicial circuits.''  In fact, 
Welsh \cite[p 180]{welsh} gives a different proof of  Lemma \ref{lm:circ-matroid} based on binary matroid theory.  We thank a referee for the references in this paragraph.

Our extension of an Eulerian trail into higher dimensions is the implementation of a simple idea: Given a regular CW-complex $K$ that obeys all pseudomanifold conditions except that its sides have positive even degrees, $K$ can be transformed into a pseudomanifold by replacing each side (of degree $r$) by $r/2$ sides of degree 2 and appropriately rearranging the attaching maps of the facets while leaving the codimension-2 skeleton unchanged.
When $K$ is 2-dimensional, a relatively simple argument suffices \cite{hk-math-mag}, but for dimensions 3 and above, we needed Lemma \ref{th:bridgeless}.
While it is nice to find such a pseudomanifold, with its natural map to $K$,  the nature of pseudomanifolds in higher dimensions \cite{pseudo} is mostly unknown.

Circlets may provide useful bases for homology of complexes such as are calculated in topological data analysis (e.g, \cite{deh, b-et-al}, cf. \cite{pk-jktr}).  Every element in ($n$-dimensional) $\mathbb{Z}_2$-homology is represented by an even ($n$-)subcomplex.  If $n=1$, then such even-degree graphs are edge-disjoint unions of cycles and so are sums of cycles, and furthermore, these cycles are {\it contained within the graph}.  This property need not hold in higher dimensions; e.g., the 2-d circlet $C(k,m)$ in Fig.~1 has no pseudosurface subcomplex.  However, every pure even subcomplex is the facet-disjoint union of a set of circlets it contains. 
Nonunique circlet-decompositions are given in \cite{hk-math-mag}.


Suppose we have a pure, even, s-connected $n$-complex $K$ such that every facet has an even number of sides and such that two facets share at most one common side.    If $\varphi: M \to K$ is an Euler cover, then $M^*$ is a simple Eulerian graph; thus, the line graph of $M^*$ is Hamiltonian  \cite[p 83]{harary}.
Hence, one can cyclically order all the sides of $M$ such that each consecutive pair of sides are in the boundary of a (unique) common facet.  
For $n \geq 3$, this tour of the sides could be prescribed by a simple closed curve in $M$ obtained by splitting the vertex of $M^*$ representing a facet $\sigma$ of $M$ into $k$ copies if the facet has $2k$ sides with all $k$ of these vertices in the interior of $\sigma$ and joining them into a simple closed non-self-intersecting curve as shown in Fig.~2.  This closed curve will pass exactly once through each side of $M$.

In the practical case of a 2-complex in 3-dimensions, something like an Euler cover appears in robotics \cite[Fig. 4]{robots}, where this concept is implicitly used to study movement trajectories separated by obstacles, via a simplicial complex based on measurements.  

\section*{Acknowledgements}

We appreciate the advice and criticisms supplied by our referees, who helped us to clarify definitions and proofs.  We are grateful for their attentive reading of the paper and for their generosity in providing specific corrections.

\end{document}